\documentclass[11pt]{article}
\usepackage[utf8]{inputenc}
\usepackage[
papersize={5.5in, 8.5in},
left=0.35in,
right=0.35in,
top=0.35in,
bottom=0.5in]{geometry}
\usepackage{graphicx}
\usepackage{amsthm}
\usepackage{rotating}
\usepackage{amsfonts,amssymb,amsmath}
\usepackage[usenames]{color}

\usepackage{array}
\usepackage{float}
\usepackage{booktabs}
\usepackage[dvipsnames,usenames,table,xcdraw]{xcolor}

\graphicspath{ {images/} }

\title{\bf 
On upper bounds for the multi-fold \\ chromatic numbers of the plane
}
\author{\bf 
\textcolor[rgb]{0,.3,1}{Jaan Parts} \\
} 
\date{\normalsize \textcolor[rgb]{0,.3,1}{Kazan, Russia, jaan\_parts@.mail.ru}}

\begin{document}

\maketitle

\pagestyle{empty}
\thispagestyle{empty}

\begin{abstract}
The multi-fold chromatic number of the plane $\chi_m$ is the smallest number of colors $k$, sufficient to color each point of the Euclidean plane in exactly $m$ colors, so that for any pair of points at a unit distance from each other, two corresponding $m$-subsets of $k$-set do not contain any common color.
We consider upper bounds for $m$-fold chromatic numbers of the plane. Our main result is that for any $m$ the inequality $\chi_m<(1+2/\sqrt3)^2\cdot m+3.501$ holds.
\end{abstract}

\section{Background}

The classical \textit{chromatic number} (CN) of a plane is introduced as follows. Each point of the Euclidean plane is associated with some color $c\in C$. A \textit{proper} coloring is one in which any pair of points located at a unit distance from each other has different colors. The smallest number of colors $k=|C|$ that provides a proper coloring is called the chromatic number of the plane $\chi$.

Even $2/3$ of a century after the definition of CN, its exact value is still unknown, and only the range of possible values of $\chi$ is known. During this time, various extensions and generalizations of the concept of CN have been invented. Now there are so many different generalizations\footnote{
There is nothing easier than coming up with a new CN: it is enough to modify the coloring rules, change the object of coloring, or introduce an additional condition for some of its characteristics.}
that it is difficult not to get confused in them. One of the examples is the \textit{fractional} CN (FCN), as well as the \textit{multi-fold} (or $m$-fold) CN (MCN). As with CN, the exact values of MCN and FCN are not known.

In the case of MCN $\chi_m$, each point of the plane is associated with a set of $m$ different colors\footnote{
One should distinguish between MCN and the \textit{list} CN (LCN). For the latter one each point is associated with its own list of allowed colors, but receives one single color in the end.}
(the value of $m$ is the same for all points), while for any pair of points at a unit distance from each other, the $m$-sets should not have common colors. Clearly, CN is a special case of MCN with $m=1$. 
FCN is defined as $\chi_f=\lim_{m\rightarrow\infty}\frac{\chi_m}{m}$. That is, MCN is intermediate between CN and FCN.

In what follows, we are only interested in the \textit{upper bounds} for MCN (best known FCN and CN upper bounds are $\chi_f<4.3599$ and $\chi\le7$). The general approach for obtaining such estimates is to create some periodic \textit{tiling}, in which the plane is densely covered $m$ times (layers) with \textit{tiles} of $k$ different colors, so that the maximum width of tiles does not exceed 1, and a distance of at least 1 is maintained between tiles of the same color. It remains only to propose the most suitable shape of the tiles.

This problem was previously considered by Grytczuk et al. \cite{gry}, where regular hexagons formed the basis of the tiling. 
We improve the known upper bounds $k$ of $\chi_m$ by considering other tile shapes.

\section{Tools}

\paragraph{Basic constructions.}
We use two basic constructions depending on the value of $m$. In both cases, the plane is colored with the same layers, shifted relative to each other. Each layer is covered without gaps by identical tiles, the centers of which form some lattice.

In the case of small $m$, each of the $m$ layers is covered by hexagonal tiles of $k$ different colors, and can be divided into identical clusters of $k$ tiles (one tile of each color). The layers are overlapped so that tiles of the same color from different layers are shifted from each other by the distance $d$ not less than 1.

In the case of large $m$, each of the $k$ layers contains tiles of only one color and can be partitioned into identical fragments, consisting of $k=k_1\times k_2$ small congruent cells, $m$ of which form one colored tile (with a width not exceeding 1), and the rest remain uncolored. 
The shape of the tiles is chosen so as to increase the proportion of the colored plane. The layers are overlapped so that each point of the plane is colored in exactly $m$ different colors.

The main task is to minimize the $k/m$ ratio for each $m$.

\paragraph{Asymptotic estimates.}
The optimal shape of the tiles depends on the number of layers allowed. For small $m$, the tiles should be angular to ensure a uniform coverage of the plane without gaps. Hexagonal tiles are best suited for this purpose. The smallest value of the $k/m$ ratio, which can be achieved theoretically by tiling the plane with regular hexagons with sides $1/2$, located at a unit distance from each other at the nodes of a hexagonal lattice, is $\tau_1=(1+2/\sqrt3)^2\approx4.642734410$. The density of the tiling depends on the angle $\alpha$ of the rotation of the hexagons relative to the lattice (or, equivalently, on the shift between the parallel sides of adjacent tiles of the same color). So the smallest $k/m$ value changes from $\tau_1$ to $16/3$ for extreme cases, where tiles are turned to each other by sides or corners, respectively.
More precisely, $k\ge\lceil\tau_1\cdot m \rceil$, due to the integer value of $k$.

At sufficiently large $m$, it becomes possible to increase the plane coloring density by using more rounded tiles \cite{pri}. So,when using circle tiles, $k/m$ can approach the threshold $\tau_2=8\sqrt3/\pi\approx4.410631163$. The best known 
bound $\tau_3=\min_{\,\theta}\left[ (1 + \cos\theta)^2/\sqrt{3}\,/\left(\pi/6-\theta+\sin(2 \theta)/2\right) \right]\approx4.359868202$ is given by Croft's construction \cite{cro} (a rounded hexagon\footnote{
In publications, one can also find such names as roundgon, or tortoise.}
with equal lengths of the straight and rounded parts of the perimeter).

One can search for $k$ in form $k=\tau_1\cdot m+\delta$ or $k=\lceil\tau_1\cdot m \rceil+\Delta$ with some real or integer correction value $\delta$ or $\Delta=\lfloor\delta\rfloor$ respectively.

\paragraph{Tile shape.}
We will focus basically on the case of small $m$, where hexagonal tiles are used. Here the problem is to find for each number of layers $m$ a tiling that minimizes the number of tiles $k$ in the repeating cluster.

To simplify the construction of the tiling,
additional restrictions are imposed on the parameters of the tiles: a) all 6 corners of the tiles are located on a circle of unit diameter; b) all distances $d$ between the nearest three tiles of the same color (from different layers) are the same.

\begin{figure}[!b]
\centering
\begin{tabular}{@{}c<{\!\!}c@{}}
 \includegraphics[scale=0.236]{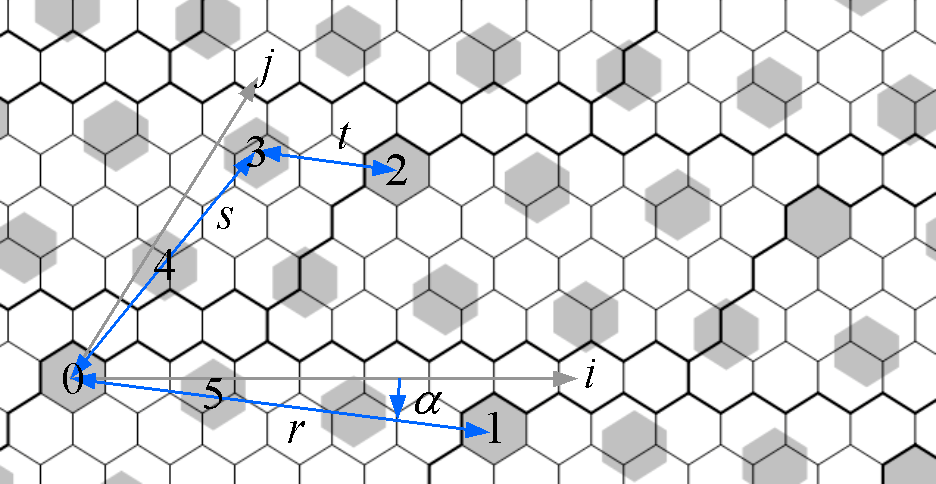} & \includegraphics[scale=0.236]{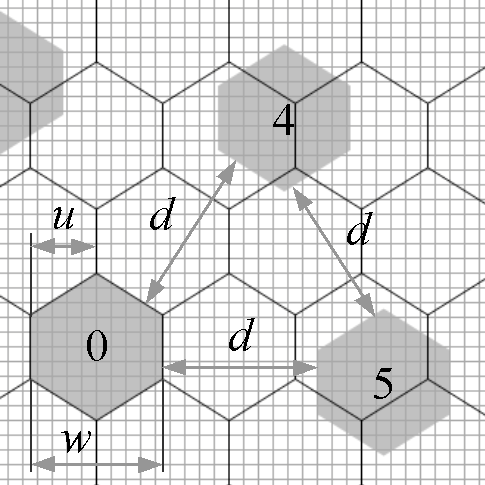} \\
\end{tabular} \par
\caption{To the definition of the tiling parameters by the example of $m=3\times2=6$ and $k=31$. The base layer tiles are shown, as well as the base color tiles (highlighted in gray). Other layers are hidden so as not to clutter up the picture, but can be unambiguously restored by the position of the base color tiles. A possible clustering option is shown in bold lines. The enlarged fragment of the tiling on the right shows the three nearest tiles of the same color on a square grid with a step 0.1.}
\label{gpar}
\end{figure}

Thereby, each tile is symmetrical about the center, has three pairs of parallel sides and is characterized by only two free parameters, for example, the width $w$ and the position of the upper corner $u$ (see Fig.\ref{gpar}). Looking ahead, we note that further complication of the tile shape (with one exception) and the introduction of additional free parameters, in our experience, did not lead to a decrease in $k$. The parallel sides of the tile are the same and form a rectangle with unit diagonals. We place one of these rectangles vertically, so that a horizontal row of tiles is formed.

\paragraph{Tiling parameters.}
Some color (tile, cluster, layer) is declared \textit{base}.

For the convenience of indexing and determining the position of the tiles, oblique coordinates $p_n=(i_n, j_n)$ are introduced, where $i_n$ and $j_n$ are integer indices of the tile $n$ along two (out of three possible) directions of the hexagonal lattice coinciding with the centers of the tiles in the base layer. The tiles of the remaining layers generally have fractional 
indices.

Adjacent tiles of the base color (from different layers) form another hexagonal lattice and are shown in gray in the figures. More specifically, $m$ adjacent tiles (one from each layer) form a parallelogram corresponding to one of the possible variants of integer decomposition $m=r\times s$.

To carry out the calculations, it is enough to consider the coordinates of six tiles of the base color, numbered from 0 to 5 (see Fig.\ref{gpar}): $p_0=(0,0)$ for the base tile; $p_1$ and $p_2$ for the base tiles of two adjacent clusters of the same layer; $p_1$ and $p_3$ for the tiles corresponding to the decomposition $m=r\times s$; 
$p_0$, $p_4$ and $p_5$ for adjacent tiles that form a triangle of the minimum perimeter; moreover, $p_5=p_1/r$, $p_4=p_3/s$, $p_3=p_2-t\cdot p_5$.

The concrete construction of the tiling is described by the parameters
characterizing the density of the tiling $\{m, k\}$, the tile shape $\{w, u\}$, the distance between the tiles $\{d\}$, the shape of the cluster $\{ i_1, j_1, i_2, j_2\}$, an arrangement of tiles of the same color $\{r, s, t\}$. 
Note that these parameters are not fully independent.

\paragraph{Computational procedure.}
We seek to minimize $k$ and, subject to that, to maximize $d$ for each $m$.
The search for the optimal tiling is performed by exhaustive search. For given $m$ the values of $k$ are sequentially increased, starting from $k=\lceil (1+2/\sqrt3)^2\cdot m \rceil$, until a tiling is found with $d\ge 1$. For each $m$, all possible variants of $\{r, s, t\}$ are checked. For each $k$, a set of possible $\{p_1, p_2\}$ is determined taking into account the relation $k=i_1\cdot j_2 - i_2\cdot j_1$. For each configuration described by integer parameters $\{m, k, i_1, j_1, i_2, j_2, r, s, t\}$, a procedure is started with variable parameters $w$ and $u$, which maximizes the value of $d$, so that $0\le u\le w\le 1\le d$. 

To reduce the search area, a number of other natural restrictions are introduced on the range of parameter values:
$j_1\le0$, $i_1\ge2\,|j_1|$ (take into account the symmetry of the tiling);
$i_1>2\,r$ (provides a horizontal gap between the tiles);
$k>2\,s\,i_1$ (sets a lower limit on the cluster size);
$0\le i_2<i_1$, $j_2\le 2\,s+1$, $i_2+j_2>2\,s$, $t<r$ (limit the position $p_2$).

\section{Tricks}

\paragraph{Extrapolation.}
If one superimpose two arbitrary tilings with parameters $\{m_1, k_1\}$ and $\{m_2, k_2\}$ one on top of the other, then one can always get a combined tiling with parameters $\{m, k\}$ using $k=k_1+k_2$ different colors and $m=m_1+m_2$ layers.

Combining the tilings allows one to obtain improved MCN estimates for some $m$ from the available MCN values for smaller $m$, as well as to extrapolate to larger $m$ values. For this, all variants of the decomposition of $m$ into terms $m=m_1+m_2$ are sorted out, and for them the minimum sum $k=k_1+k_2$ is determined. Practice shows that such extrapolation gives good results (see Section~\ref{sres}).

\paragraph{Wavy sides.}
If the parallel sides of adjacent tiles of the same color are shifted by more than half the length of a side\footnote{
Notice that if all tiles have the same shape then each side of the tile must be symmetrical about its center. Thus the position of the center relative to the ends of the side is fixed.},
one can use the trick of 
drawing the relevant part of this side as an arc with a radius of 1 centered at
the nearest corner of the adjacent tile to slightly reduce the distance between the tiles. The largest effect (reduction of distance by about $0.8\%$) is observed when the sides of adjacent tiles are shifted by $3/4$ of their length.

Potentially, the effect of wavy sides can work independently for any pair of parallel sides.
Since, with increasing $m$, there is a general tendency towards a decrease in the angle $\alpha$ (tiles tend to take the shape of a regular hexagon and turn sides to each other), one can expect the appearance of wavy sides only for small $m$ (practically for $m<40$).

\section{Results}\label{sres}

\paragraph{Small $m$.}
Let's start with the first construction. Here for a particular pair $\{m, k\}$, it is sufficient to get the only tiling with $d\ge 1$. Usually, one can find several different tilings with this property\footnote{
Also note that the chosen notation system has redundancy, and the same tiling can be specified in different ways.}. 
We show those that maximize $d$. Along with $k$, we use the deviation $\delta$ from the theoretical minimum $\tau_1$ for hexagonal tiles.

We did a full search for all $m\le 100$, as well as for $m\le 2000$ with some simplifications and gaps. Further search is devoid of practical sense, bearing in mind the increase in its labor intensity, the possibility of simple extrapolation, as well as the transition to the second construction.

Most of the obtained values $\delta=k-\tau_1$ are within the range $\delta\in(0.35, 3)$. The exceptions are $m\in\{7,48,132,135\}$. For all $m>135$ we got $\delta<2.74$. The largest value $\delta<3.501$ is observed at $m=7$.

The lower threshold $\delta>0.35$ is explained by the fact that, if the shape of the tiles approaches a regular hexagon, a row of tiles of the same color ($r$-axis on Fig.\ref{gpar}) should have a slight slope with respect to a row of tiles of the same layer ($i$-axis). This means that the value $\tau_1$ is unattainable in this construction for finite $m$.

Oddly enough, we managed to find the only (but practically important) case when the use of wavy sides leads to a decrease in $k$: for $m=6$, it is possible to reduce $k$ from 31 to 29 (due to reduction of $d$ by only $0.15\%$).

\begin{figure}[H]
\centering
\begin{tabular}{@{}c<{\!\!}c<{\!\!}c@{}}
7/1 & \includegraphics[scale=0.29]{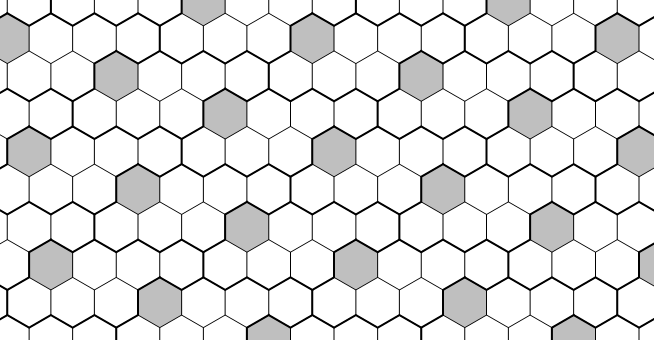} & \includegraphics[scale=0.29]{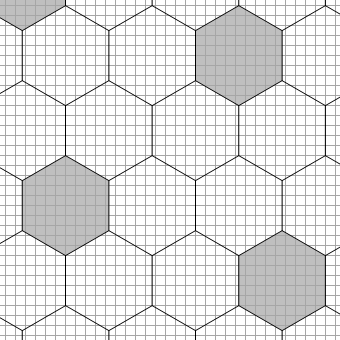} \\ [4pt]
11/2 & \includegraphics[scale=0.29]{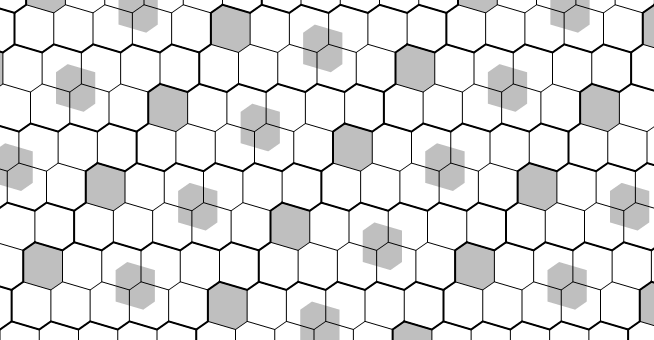} & \includegraphics[scale=0.29]{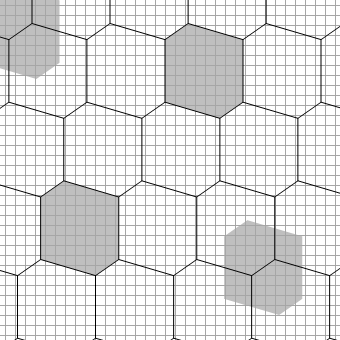} \\ [4pt]
16/3 & \includegraphics[scale=0.29]{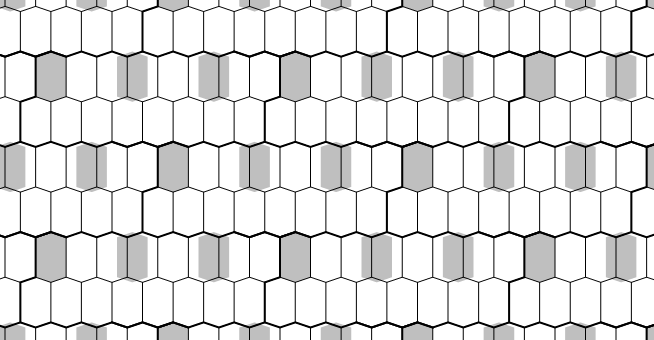} & \includegraphics[scale=0.29]{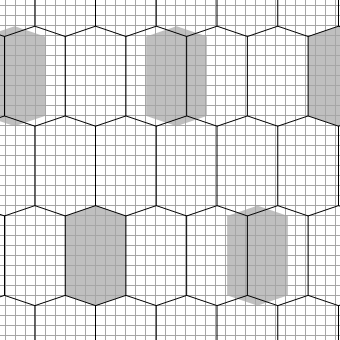} \\ [4pt]
21/4 & \includegraphics[scale=0.29]{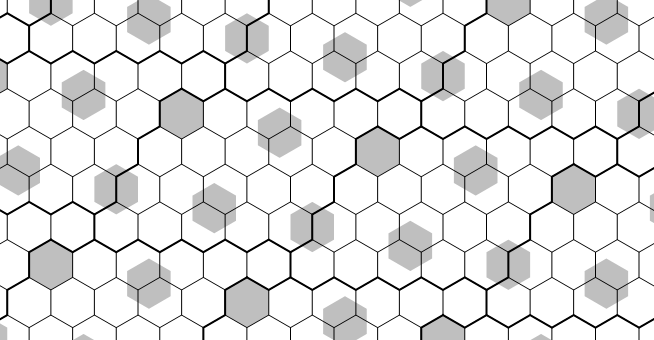} & \includegraphics[scale=0.29]{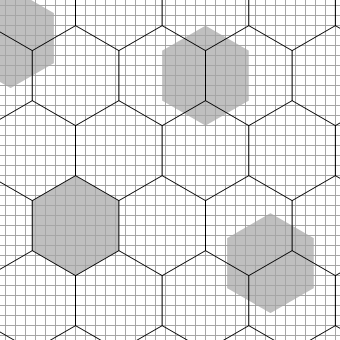} \\ [4pt]
26/5 & \includegraphics[scale=0.29]{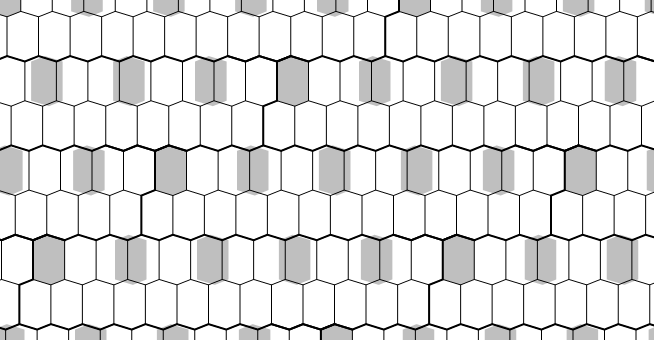} & \includegraphics[scale=0.29]{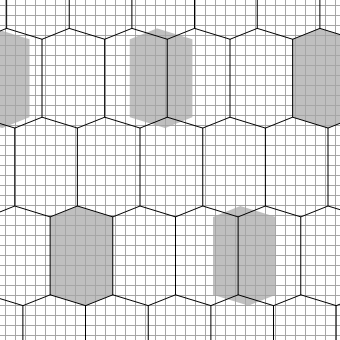}
\end{tabular} \par
\caption{Tilings for $m=1\ldots5$ and best found values of $k/m$.}
\label{glo1}
\end{figure}

\begin{figure}[H]
\centering
\begin{tabular}{@{}c<{\!\!}c<{\!\!}c@{}}
29/6 & \includegraphics[scale=0.29]{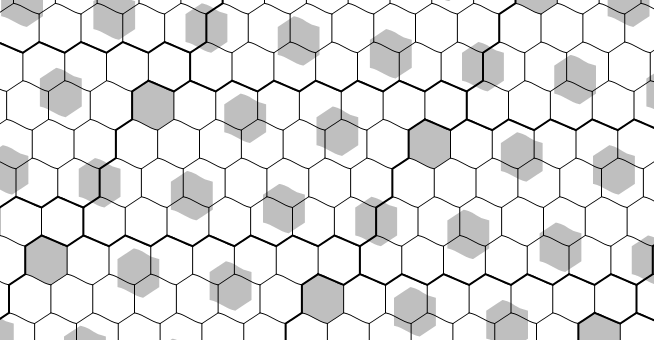} & \includegraphics[scale=0.29]{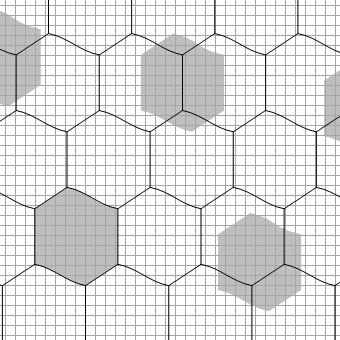} \\ [4pt]
36/7 & \includegraphics[scale=0.29]{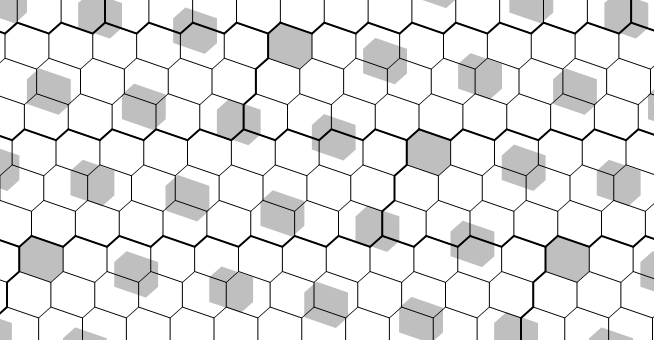} & \includegraphics[scale=0.29]{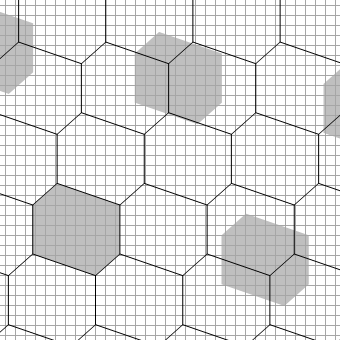} \\ [4pt]
39/8 & \includegraphics[scale=0.29]{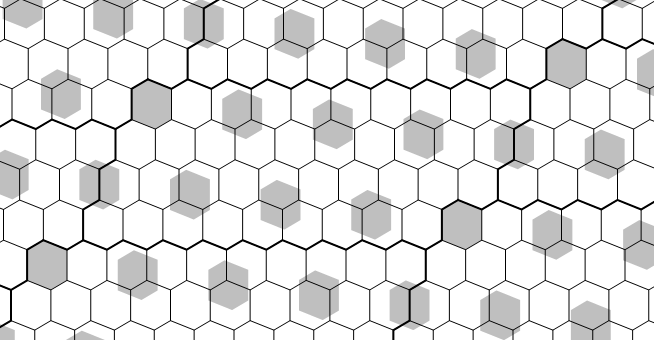} & \includegraphics[scale=0.29]{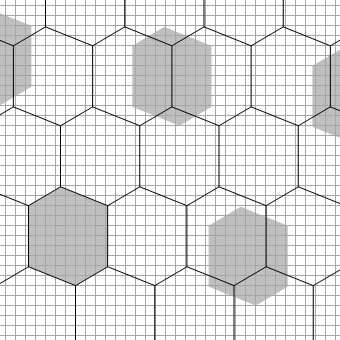} \\ [4pt]
43/9 & \includegraphics[scale=0.29]{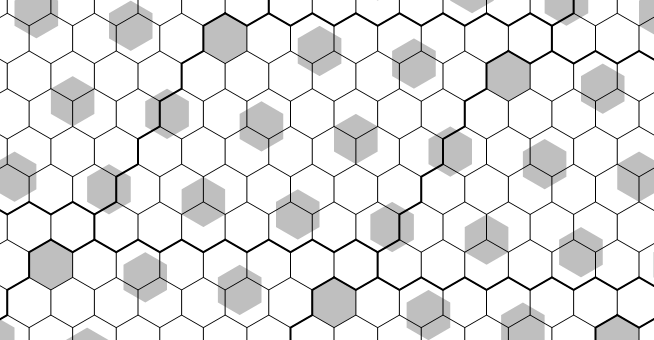} & \includegraphics[scale=0.29]{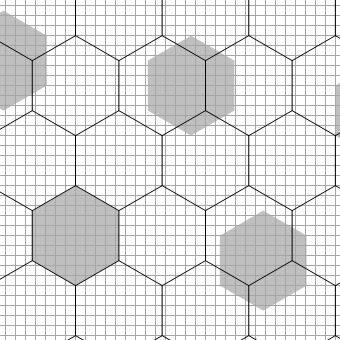} \\ [4pt]
49/10 & \includegraphics[scale=0.29]{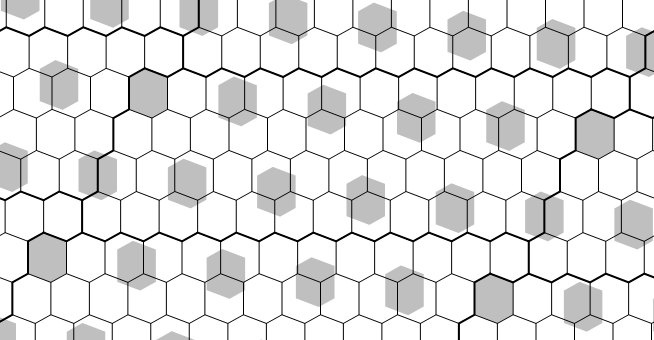} & \includegraphics[scale=0.29]{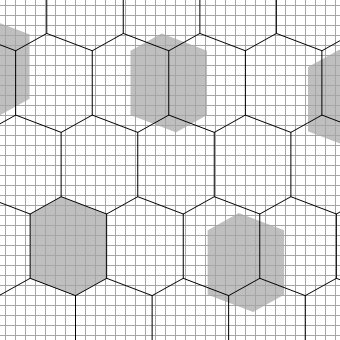} 
\end{tabular} \par
\caption{Tilings for $m=6\ldots10$ and best found values of $k/m$.}
\label{glo2}
\end{figure}

Figures \ref{glo1} and \ref{glo2} show some examples of tilings that give the best found upper bounds of MCN for $m\le 10$. If one look closely, one can see the wavy sides in the case $k/m=29/6$. The discovery of this tiling increases the efficiency of the extrapolation trick. In the range $100<m\le 2000$, we got 1456 explicit tilings with $\delta<3$. The remaining 444 estimates of $\delta$ are obtained by extrapolation. Besides, extrapolation allowed to decrease $\delta$ by 1 and 2 for 87 and 12 cases of explicit tiling, respectively.

\newpage

\begin{table}[H]
{
\caption{Characteristics of the best found tilings.}
\label{tbest}
\smallskip

{
\centering
\small
\begin{tabular}{@{}>{\;}r<{\;}@{}| *{2}{@{}>{\;}r<{\;}@{}}*{2}{@{}>{\;}c<{\;}@{}}| *{3}{@{}>{\;}r<{\;}@{}}@{\,}r<{\;}@{}| @{}>{\;}r<{\;}@{}>{\;\,}r<{\;}@{}r<{\;}@{}|*{3}{@{}>{\;}c<{\;}@{}}}
\hline
$m$&$k$&$k^*$&   $k/m$& $\delta$&$i_1$&$j_1$&$i_2$&$j_2$&$r$&$s$&$t$& $w$&   $u$&      $d$\\
\hline
\hline
 1&   7&   7& 7.0000& 2.3573&  3& -1&  1&  2&  1& 1&  0& 0.866025& 0.433013& 1.3229\\
 2&  11&  12& 5.5000& 1.7145&  5& -1&  1&  2&  2& 1&  0& 0.780412& 0.231291& 1.0550\\
 3&  16&  16& 5.3333& 2.0718&  8&  0&  3&  2&  3& 1&  1& 0.607271& 0.303636& 1.0121\\
 4&  21&  24& 5.2500& 2.4291&  5& -1&  1&  4&  2& 2&  0& 0.866025& 0.433013& 1.0825\\
 5&  26&  32& 5.2000& 2.7863& 13&  0&  3&  2&  5& 1&  1& 0.626508& 0.272476& 1.0024\\
 6&  29&  32& 4.8333& 1.1436&  7& -1&  1&  4&  3& 2&  0& 0.831425& 0.322944& 1.0009\\
 6&  29&  32& 4.8333& 1.1436&  7& -1&  1&  4&  3& 2&  0& 0.831379&  0.327403& 0.9994\\
 6&  31&  32& 5.1667& 3.1436&  7& -1&  3&  4&  3& 2&  1& 0.894452&  0.362876& 1.0716\\
 7&  36&  37& 5.1429& 3.5009& 16& -3& 12&  0&  7& 1&  5& 0.871616& 0.243410& 1.0163\\
 8&  39&  45& 4.8750& 1.8581& 19& -2& 10&  1&  8& 1&  4& 0.792258& 0.321035& 1.0091\\
 9&  43&  49& 4.7778& 1.2154&  7& -1&  1&  6&  3& 3&  0& 0.866025& 0.433013& 1.0104\\
10&  49&  55& 4.9000& 2.5727& 12& -1&  1&  4&  5& 2&  0& 0.765621& 0.310729& 1.0097\\
11&  53&  72& 4.8182& 1.9299& 25& -4&  7&  1& 11& 1&  3& 0.928440& 0.477605& 1.0080\\
12&  57&  63& 4.7500& 1.2872&  9& -1&  3&  6&  4& 3&  1& 0.885280& 0.397978& 1.0071\\
13&  63&  87& 4.8462& 2.6445& 30& -3& 21&  0& 13& 1&  9& 0.866025& 0.433013& 1.0326\\
14&  67&  80& 4.7857& 2.0017& 31& -4&  9&  1& 14& 1&  4& 0.931056& 0.445051& 1.0034\\
15&  71&  77& 4.7333& 1.3590& 11& -1&  5&  6&  5& 3&  2& 0.898147& 0.379766& 1.0018\\
16&  77&  81& 4.8125& 2.7162&  9& -1&  5&  8&  4& 4&  2& 0.911429& 0.475233& 1.0205\\
17&  81& 111& 4.7647& 2.0735& 39& -4& 30& -1& 17& 1& 13& 0.872477& 0.475928& 1.0171\\
18&  86&  91& 4.7778& 2.4308& 14& -2&  1&  6&  6& 3&  0& 0.866025& 0.433013& 1.0104\\
19&  91& 123& 4.7895& 2.7880& 43& -3& 16&  1& 19& 1&  7& 0.866025& 0.433013& 1.0256\\
20&  95&  99& 4.7500& 2.1453& 22& -1&  7&  4& 10& 2&  3& 0.868899& 0.341116& 1.0086\\
21&  99& 112& 4.7143& 1.5026& 48& -5& 39& -2& 21& 1& 17& 0.875384& 0.503138& 1.0057\\
22& 103& 120& 4.6818& 0.8598& 49& -3& 18&  1& 22& 1&  8& 0.861306& 0.398575& 1.0027\\
23& 109& 150& 4.7391& 2.2171& 51& -2& 29&  1& 23& 1& 13& 0.858033& 0.367470& 1.0126\\
24& 113& 117& 4.7083& 1.5744& 53& -3& 20&  1& 24& 1&  9& 0.878789& 0.419435& 1.0094\\
25& 119& 121& 4.7600& 2.9316& 11& -1&  9& 10&  5& 5&  4& 0.922862& 0.495630& 1.0083\\
26& 123& 144& 4.7308& 2.2889& 29& -2& 18&  3& 13& 2&  8& 0.883319& 0.478525& 1.0135\\
27& 127& 140& 4.7037& 1.6462& 20& -1&  7&  6&  9& 3&  3& 0.866025& 0.433013& 1.0104\\
28& 131& 144& 4.6786& 1.0034& 61& -2& 35&  1& 28& 1& 16& 0.871101& 0.367408& 1.0004\\
29& 137& 189& 4.7241& 2.3607& 65& -4& 18&  1& 29& 1&  8& 0.869685& 0.484585& 1.0128\\
30& 141& 143& 4.7000& 1.7180& 33& -1&  9&  4& 15& 2&  4& 0.864112& 0.407710& 1.0098\\
31& 145& 201& 4.6774& 1.0752& 69& -5& 29&  0& 31& 1& 13& 0.882489& 0.495857& 1.0018\\
32& 151& 162& 4.7188& 2.4325& 70& -1& 11&  2& 32& 1&  5& 0.864953& 0.409509& 1.0143\\
33& 155& 168& 4.6970& 1.7898& 72& -1& 11&  2& 33& 1&  5& 0.863131& 0.387042& 1.0083\\
34& 159& 184& 4.6765& 1.1470& 75& -2& 42&  1& 34& 1& 19& 0.854345& 0.418193& 1.0056\\
35& 165& 176& 4.7143& 2.5043& 76& -1& 13&  2& 35& 1&  6& 0.876328& 0.469624& 1.0131\\
36& 169& 169& 4.6944& 1.8616& 13&  0&  0& 13&  6& 6&  0& 0.866025& 0.433013& 1.0104\\
37& 173& 240& 4.6757& 1.2188& 80& -1& 13&  2& 37& 1&  6& 0.875603& 0.426892& 1.0061\\
38& 177& 205& 4.6579& 0.5761& 82& -1& 13&  2& 38& 1&  6& 0.874272& 0.406629& 1.0016\\
39& 183& 203& 4.6923& 1.9334& 83& -1& 17&  2& 39& 1&  8& 0.901827& 0.448498& 1.0059\\
40& 187& 198& 4.6750& 1.2906& 85& -1& 17&  2& 40& 1&  8& 0.900383& 0.428203& 1.0022\\
\hline
\end{tabular}

}
}
\end{table}

\newpage

The best results for $m\le 40$ are listed in Table \ref{tbest}. Here, for each $m$, the best found values of $k$ and $d$ are shown, as well as the corresponding tiling parameters $\{i_1, j_1, i_2, j_2, r, s, t, w, u\}$. For comparison, the estimates of $k^*$ are given according to the data from \cite{gry}. They are not only weaker, but also do not provide convergence to some threshold:
for arbitrarily large $m$ there are estimates of $k^*/m$ here both approaching $\tau_1$ (if $m$ has large divisors) and exceeding $3+2\sqrt3\approx6.46$ (if $m$ is prime).
For $m=6$, three solutions are shown: one with wavy sides and two with straight sides with $k=29$ and $k=31$ (option $k=30$ does not give valid solutions even with wavy sides). The second solution with $k=29$ is not valid, since it does not satisfy the condition $d\ge 1$, and is presented for comparison only.

\begin{figure}[!b]
\centering
\begin{tabular}{@{}c@{}}
\includegraphics[scale=0.346]{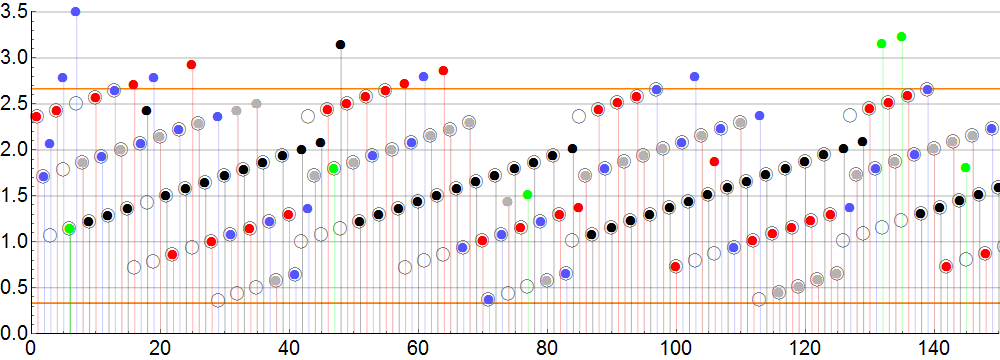}  \\ [5pt]
\includegraphics[scale=0.346]{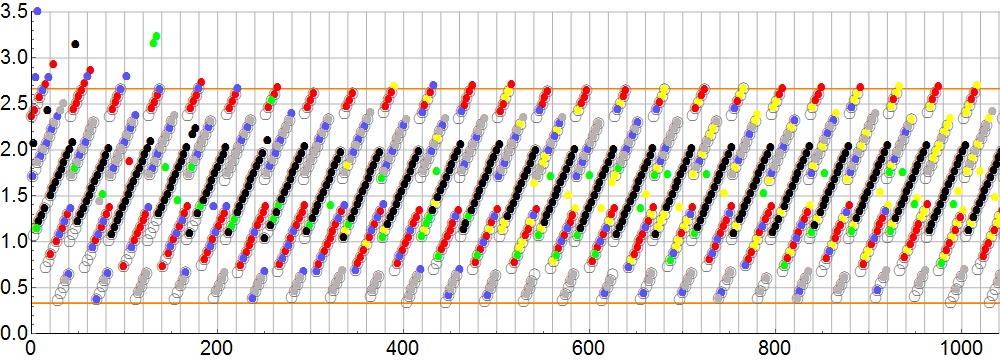} \\
\end{tabular} \par
\caption{The best found values $\delta=k-\tau_1\cdot m$ for different $m$ (for two different $m$ scales). Points show the actual values of $\delta$, circles – their probable estimates $\hat{\delta}$ (see Section~\ref{spuz}). The color indicates that the result belongs to a certain class. The sequences $m\bmod 3=\{0,1,2\}$ are shown in black, red and gray, respectively. Prime $m$ are highlighted in blue. Yellow and green show the values obtained and improved, respectively, by tricks (wavy sides and extrapolation). Colors are listed in ascending order of priority. The orange lines correspond to $\delta=\{1/3,8/3\}$.}
\label{gmed}
\end{figure}

\paragraph{Medium $m$.}
The results for $m\le 150$ and $m\le 1040$ are shown in Fig.\ref{gmed} in the form of dependence $\delta(m)$. Here the periodicity in $m$ is clearly observed. One can notice periods of about 3 and 41.75. 
The both are explained by the fact that the threshold $\tau_1$ is close to $14/3$, so that a) an increase in $m$ by 3 does little to change the fractional part $f$ of the $k/m$ ratio, b) $f$ has period $1/(14/3-\tau_1)\approx41.78461$.

With further extrapolation of more or less known $k$ for all $m\le2000$, it is possible to obtain estimates $\Delta=\lfloor\delta\rfloor\le 2$ for $m<5019$ and $\Delta\le 3$ for $m<11829$.
If we restrict the known range to $m\le500$, then we get $\Delta\le 2$ for $m<1425$ and $\Delta\le 3$ for $m<2704$. 
If we initially take $m\le40$, then after extrapolation $\Delta\le 4$ for $m<101$ and $\Delta\le 15$ for $m<893$. In Fig.\ref{gmed}, the values obtained or improved by extrapolation are highlighted in color.

\begin{figure}[!b]
\centering
\includegraphics[scale=0.35]{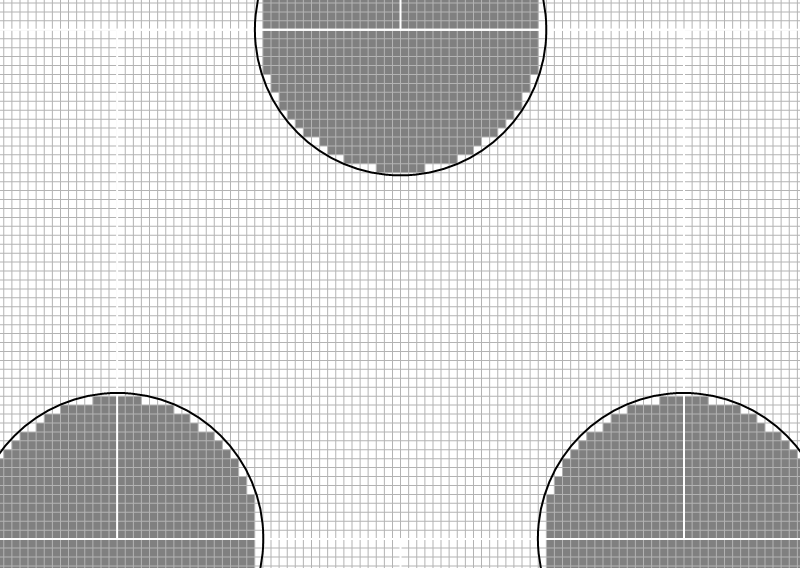}
\caption{The tiling with $m=860$ and $k/m<\tau_1$.}
\label{glar}
\end{figure}

\paragraph{Large $m$.}
In the case of large $m$ (about 1000 or more), the hexagonal tiles become ineffective.

For any given threshold value $\tau>\tau_3$, one can find $m$ for which the tiling of the second type will certainly provide $k/m<\tau$. However, such estimate is rather rough, and for $\tau=\tau_1$ gives a value $m$ of the order $10^4$. Alternatively, one can take the densest coverage of the plane with circles or Croft's smoothed hexagons, and consider different ways of filling these tiles with small rectangles (or hexagons), choosing the option with minimal loss of area around the tile perimeter.

So, for $m=860$, we found a tiling with $k=70\times 57$ (the corresponding sides of the rectangular cells are $0.0277851\times0.0306606$) and $k/m=3990/860\approx 4.63953$ (see Fig.\ref{glar}), for $m=992$ we obtained a tiling with $k=74\times 62$ and $k/m=4588/992\approx 4.625$. 
In other words, for $m=860$ we get $\delta<0$. Combining now the constructions of the first and second types\footnote{
Note that a construction of the second type exists for every $m$, since, if necessary, one can leave several rectangular cells uncolored.}
and using the trick with extrapolation, it is easy to make sure that for all subsequent $m$ the relation $k/m<\tau_1$ also holds, which means that the previously obtained estimate $\delta\le 3.501$ is valid for all $m$.

\section{Puzzles}\label{spuz}

If we take a closer look at Fig.\ref{gmed}, we will notice the following curious thing. Most of the values of  $\delta$ fit into the range $\delta\in[0.35, 2.75]$, moreover, if we divide all $m$ into three equivalence classes ($m\bmod 3$), then most of $\delta$ fall into the following ranges:
$$\delta\in\begin{cases}
[1.05, 2.05], & m\equiv 0\pmod 3\\
[0.75, 1.4]\cup[2.4, 2.75], & m\equiv 1\pmod 3\\
[0.35, 0.7]\cup[1.7, 2.35], & m\equiv 2\pmod 3
\end{cases}$$

If we round off these values to $1/3$, we get the following rough empirical formula for the probable value of $\delta$:
$$\hat{\delta}=f+\begin{cases}1, &m\equiv0\pmod3\\ \lfloor3f+m\rfloor\bmod3, &m\not\equiv0\pmod3\end{cases}$$

Here $f=\delta-\lfloor\delta\rfloor=\lceil\tau_1\cdot m\rceil-\tau_1\cdot m$ is a fraction part of $\delta$.

For $m\le 2000$, this rule is violated in about $13.5\%$ of cases. 
If we discard the initial and the rough search zones, then in a narrower range $m\in[150, 900]$ the rule is violated in $10.5\%$ of cases.
Interestingly, there are not only positive, but also negative deviations from the estimate $\hat{\delta}$ (about $3.3\%$ of cases). 
The question arises: what explains this rule? First of all, what causes the "unit holes" in $\delta$ for $m\not\equiv 0 \pmod 3$?
Maybe for such combinations of $k$ and $m$, the tiles are lined up so badly that for any options of the tiling parameters, they do not allow to obtain a valid solution. However, such an explanation explains little, and only states a fact. Something better is required here.

Among the riddles is the only successful case of a trick with wavy sides.
One can try to find other ways to optimize the tiling. For example, is less regular tiling beneficial with tiles of different shapes in the same layer? Is it possible to find an intermediate construction for medium $m$, combining an $m$-layer structure (as for hexagonal tiles) and the tiles of a more complex shape (like smoothed hexagons), and to cross the threshold $\tau_1$ much sooner, say, closer to $m\approx100$?

One can notice that the first construction for small $m$ gives convergence to the threshold $\tau_1$ with the rate $1/m$. The second construction for large $m$ gives a convergence to the smaller threshold $\tau_3$, but with slower rate $1/\sqrt{m}$. Is it possible to improve the convergence in the second case as well?

\vspace{2mm}

The article is finished. But I don't like to leave the empty spaces. Therefore, a short story.

Initially, when I had only a few hand-made tilings in my hands, I conceived an ideal mathematical article with a minimum of text, or maybe with no text at all to break the record set by Conway and Soifer in their monumental work \cite{con}. In particular, it would immediately relieve me of the need to translate it. But then, as the material began to accumulate, because my computer also wanted to take part in the creative process, I had to revise my initial views, and add some clarifications. Hope you found them interesting.
This was facilitated by A. de Grey and A. Soifer.

\end{document}